\documentclass[a4paper, 12pt]{amsart}
\usepackage[utf8]{inputenc}

\usepackage{fourier}

\usepackage{color}
\usepackage[colorlinks=true,linkcolor=colorref,citecolor=colorcita,urlcolor=colorweb]{hyperref}
\definecolor{colorcita}{RGB}{21,86,130}
\definecolor{colorref}{RGB}{5,10,177}
\definecolor{colorweb}{RGB}{177,6,38}

\usepackage{amsmath}
\usepackage{amsthm}
\usepackage{xcolor}
\usepackage{soul}

\newtheorem{thm}{Theorem}[section]
\newtheorem{cor}[thm]{Corollary}
\newtheorem{lem}[thm]{Lemma}
\newtheorem{prop}[thm]{Proposition}
\newtheorem{exa}[thm]{Example}
\newtheorem{rem}[thm]{Remark}

\DeclareFontFamily{U}{mathx}{\hyphenchar\font45}
\DeclareFontShape{U}{mathx}{m}{n}{
      <5> <6> <7> <8> <9> <10>
      <10.95> <12> <14.4> <17.28> <20.74> <24.88>
      mathx10
      }{}
\DeclareSymbolFont{mathx}{U}{mathx}{m}{n}

\DeclareMathAccent{\widecheck}{0}{mathx}{"71}

\newcommand{\cd}[1]{{{#1}^{iv}}}
\newcommand{\qd}[1]{{{#1}^{v}}}
\newcommand{\sd}[1]{{{#1}^{vi}}}

\def\M{\mathcal {M}}
\def\H{\mathcal {H}}

\DeclareMathOperator*{\wslim}{\it{w}*-\lim}

\newcommand{\zC}{\mathbb{C}}
\newcommand{\zN}{\mathbb{N}}

\begin{document}

\title[Homomorphisms on non-symmetrically regular spaces]{Homomorphisms on algebras of analytic functions on non-symmetrically regular spaces}
\author{Daniel Carando}
\address{Departamento de Matem\'{a}tica, Facultad de Ciencias Exactas y Naturales, Universidad de Buenos Aires, (1428) Buenos Aires,
Argentina and CONICET} \email{dcarando@dm.uba.ar}

\author{Ver\'onica Dimant}
\address{Departamento de Matem\'{a}tica y Ciencias, Universidad de San
Andr\'{e}s, Vito Dumas 284, (B1644BID) Victoria, Buenos Aires,
Argentina and CONICET} \email{vero@udesa.edu.ar}

\author[Tom\'as Rodr\'{i}guez]{Jorge Tom\'as Rodr\'{i}guez}
\address{Departamento de Matem\'{a}tica and NUCOMPA, Facultad de Cs. Exactas, Universidad Nacional del Centro de la Provincia de Buenos Aires, (7000) Tandil, Argentina and CONICET}
\email{jtrodrig@dm.uba.ar}

\begin{abstract}
    We study homomorphisms on the algebra of analytic functions of bounded type on a Banach space. When the domain space lacks symmetric regularity, we show that in every fiber of the spectrum there are evaluations (in higher duals) which do not coincide with evaluations in the second dual. We also consider the commutativity of convolutions between evaluations. We show that in some Banach spaces $X$ (for example, $X=\ell_1$) the only evaluations that commute with every other evaluation in $X''$ are those in $X$. Finally, we establish conditions ensuring the symmetry of the canonical extension of a symmetric multilinear operator (on a non-symmetrically regular space) and present some applications.
\end{abstract}

\thanks{This work was partially supported by CONICET PIP 11220130100329CO, 
CONICET PIP 11220200101609CO,   
ANPCyT PICT 2018-04104, 
and ANPCyT PICT 2018-04250. 
}
\subjclass[20120]{46G20, 46J15, 47H60.}

\maketitle

The study of the spectrum of the Fr\'{e}chet algebra  $\H_b(X)$ of entire functions of bounded type began with the renowned article \cite{AroColGam91} by Aron, Cole and Gamelin. In that work, the hypothesis of $X$ being symmetrically regular (see the definition below) was identified as relevant for the description of its spectrum. From there on, all the subsequent research on this topic (see, for instance,  \cite{AroGalGarMae96,CarGarMae05,CarGarMaeSev2012,CarMur,Din99}) has been mainly devoted to the symmetrically regular case.

Our aim here is to follow the unexploited path and study the behavior of homomorphisms in $\M_b(X)$ for $X$ a non-symmetrically regular Banach space.

The spectrum of $\H_b(X)$ is the set
\[
\M_b(X)=\{\varphi:\H_b(X)\to\mathbb C,\ \varphi \textrm{ is a nonzero algebra homomorphism}\}.
\]
For each  $x\in X$, we can define the \emph{evaluation homomorphism at}
 $x$ as $\delta_x(f)=f(x)$, which clearly belongs to $\M_b(X)$.
By \cite{AroBer78}, there is a canonical extension $[f\mapsto \tilde f]$  from $\H_b(X)$ to $\H_b(X'')$, the so-called  Aron-Berner extension (see definition below). Through this extension we can also define \emph{evaluations} at elements of the bidual:  for  $x''\in X''$ we set 
 \begin{equation}\label{eval en bidual}\tilde{\delta}_{x''}(f):=\tilde{f}(x'').\end{equation}
Since this extension is continuous and multiplicative we have $\tilde{\delta}_{x''}\in\M_b(X)$ for each
$x''\in X''$. This allows us to consider $X''$ as a subset of $\M_b(X)$. There is also a natural projection from $\M_b(X)$ to $X''$: as $X'$ is contained in $\H_b(X)$, the mapping $\pi:\M_b(X)\to X''$ is defined by 
\begin{equation}\label{defi-pi}
\pi(\varphi)=\varphi_{\big|_{X'}}, \end{equation}
i.e.,  $\pi(\varphi)$ is the restriction of the
homomorphism $\varphi$ to $X'$. For  $x'' \in X''$,  {\it the fiber} of the spectrum of $\M_b(X)$
over  $x''$ is the set
 $$\mathcal{M}_{x''}=\{\varphi \in \M_b(X),\   \pi(\varphi)=x''\}=\pi^{-1}(x'').$$
 It is easy to check that $\pi(\tilde\delta_{x''})=x''$ or, in other words, that  $\tilde\delta_{x''}$ belongs to $\mathcal{M}_{x''}$. 
  We thus have the following commutative diagram:
\vspace{0.5cm}
\begin{center}
	\begin{picture}(130,90)
	\put(0,85){$X''\,\,$ } \put(112,85){$\M_b(X)$} \put(120,15){$X''$}
	\put(18,88){\vector(1,0){95}} \qbezier(18,88)(15,88)(15,91) \qbezier(15,91)(15,94)(18,94)
	\put(123,81){\vector(0,-1){55}} \put(10,80){\vector(2,-1){110}}
	\put(62,90){$\tilde\delta$}
	\put(126,50){$\pi$}\put(55,45){$id$}
	\end{picture}
\end{center}

 The question on the existence of elements in $\mathcal{M}_{x''}$ other than $\tilde\delta_{x''}$ has been widely considered. In particular, one may wonder if we can get new homomorphisms by repeating the extension and evaluation procedure. Namely, let $\tilde{\tilde f}$ be the extension of $\tilde f$ to the fourth dual $X^{iv}$ and take $x^{iv}\in X^{iv}$. We can define $$\tilde{\tilde{\delta}}_{x^{iv}}(f)=\tilde{\tilde{f}}(x^{iv}).$$ 
 This is clearly an element of the spectrum, but it is not clear that it is a new one. It may happen that there exists some $x''\in X''$ such that \begin{equation}\label{eq-noesnuevo}
     \tilde{\tilde{\delta}}_{x^{iv}}(f)=\tilde{\delta}_{x''}(f)
 \end{equation} for all $f\in \H_b(X)$. In fact, in \cite{AroGalGarMae96} it is shown that this is always the case when $X$ is symmetrically regular. 
 Note that if \eqref{eq-noesnuevo} holds then  one must have $x''=\pi(x^{iv})$.
Also in \cite{AroGalGarMae96}, the authors show that in every non-symmetrically regular space there exist some
$x^{iv}\in X^{iv}$ for which $\tilde{\tilde{\delta}}_{x^{iv}}(f)\ne\tilde{\delta}_{\pi(\cd{x})}(f)$ for certain $f\in \H_b(X)$. 
In Section \ref{Sect:Cuarto dual} we find conditions on $X$ ensuring that on every fiber there are  new evaluations at points of the fourth dual; that is, for \textbf{every} $x''\in X''$  there is $\cd{x}\in \cd{X}$ such that 
$$\tilde{\tilde{\delta}}_\cd{x}\in \mathcal{M}_{x''}\,\,\,\text{and}\,\,\,\tilde{\tilde{\delta}}_{\cd{x}}\neq \tilde{\delta}_{x''}.$$ These conditions are fulfilled by $\ell_1$ and,  more generally, for any $X$ having a complemented copy of $\ell_1$. For general non-symmetrically regular Banach spaces we show that  for  \textbf{every} $x''\in X''$  there is $\sd{x}\in \sd{X}$ such that 
$\tilde{\tilde{\delta}}_\sd{x}\in \mathcal{M}_{x''} $ and $\tilde{\tilde{\delta}}_{\cd{x}}\neq \tilde{\delta}_{x''}.$  In other words, in every fiber there are  homomorphisms given by evaluations in the sixth dual which do not coincide with evaluations in the bidual.
We also address a question on approximation of homomorphisms by evaluations. From \cite{DavGam89} we see that evaluations in the bidual $X''$ are limits of evaluations in $X$. 
We observe in Proposition \ref{prop-corona} that the same happens with evaluations at points of $\cd{X}$ and of any higher dual: they can be approximated by evaluations in $X$ and therefore they do not belong to the so-called corona of the spectrum.

In \cite{AroColGam91}  the authors also defined and studied the convolution operation $\ast$ (see definition below), which endows the spectrum with a semigroup structure. 
 They presented some useful formulas to compute this convolution and showed that  this operation makes the spectrum  a non-commutative semigroup.  Moreover, they give equivalent formulations for commutativity between evaluations in the bidual, i.e., conditions equivalent to $$\tilde{\delta}_{z_1}\ast \tilde{\delta}_{z_2}=\tilde{\delta}_{z_2}\ast\tilde{\delta}_{z_1}$$ for $z_1,z_2\in X''$.
 In Section \ref{Sect:convolucion} we continue the study of the convolution operation on non-symmetrically regular Banach spaces. We show that if $X$ is $\ell_1$ or $F\oplus F'$, with $F$ non-reflexive, a homomorphism $\tilde{\delta}_{z}$ commutes with every other evaluation if and only if $z\in X$.

 A careful reading of \cite{AroGalGarMae96} shows that a key ingredient used by the authors to find new homomorphisms given by evaluations at  points of $\cd{X}$ is  the existence of a symmetric multilinear map whose canonical extension is not symmetric. This is the reason why on symmetrically regular Banach spaces the points of $\cd{X}$ do not yield to  new homomorphisms: a space is symmetrically regular if and only if the canonical extension of every symmetric multilinear map is again symmetric. In Section \ref{Sect:AB ext} we investigate when the canonical extension of a particular symmetric multilinear map is symmetric. On Theorem \ref{th AB sim} we prove  several characterizations of this fact which lead to some applications. In particular, we show that if the canonical extension of a  symmetric multilinear form is symmetric, then so are all subsequent extensions to higher duals.

\subsection*{Definitions and notation}

 For a Banach space $X$ we write $B_X$ for the open unit  ball and $i_X:X\to X''$ for the canonical embedding of $X$ into its bidual. Note that, in this case, the transpose  $i_X': X'''\to X' $ is just the restriction operator, which can be seen as a natural projection from $X'''$ onto its subspace $i_{X'}(X')$.  It is important to recall that, in general, the bitranspose $i_X'':X''\to X^{iv}$ is not equal to the canonical embedding of the bidual into the fourth dual $i_{X''}:X''\to X^{iv}.$

Given Banach spaces $X$ and $E$  we say that a function $P:X\to E$ is a {\em continuous $m$-homogeneous polynomial} if there exists a unique continuous symmetric $m$-linear mapping $\widecheck{P}$ such that $P(x) = \widecheck{P}(x,\dots, x)$. We write $\mathcal{P}(^mX,E)$, $\mathcal{L}(^mX,E)$ and $\mathcal{L}_s(^mX,E)$ for spaces of  $m-$homogeneous polynomials, $m-$linear mappings and symmetric $m$-linear mappings from $X$ to $E$, respectively. When $E=\mathbb C$ we just write  $\mathcal{P}(^mX), \mathcal{L}(^mX)$ and $\mathcal{L}_s(^mX)$ for these spaces.

A mapping $f:X\to \zC$ is said to be {\em holomorphic} if for every $x_0\in X$ there exists a sequence $(P_mf(x_0))$, with each $P_mf(x_0)$ a continuous $m$-homogeneous polynomial, such that the series
\[
f(x)=\sum_{m=0}^\infty P_mf(x_0)(x-x_0)
\] converges uniformly in some neighborhood of $x_0$.

A holomorphic function $f:X\to \mathbb C$ is said to be {\em of bounded type} if it maps bounded subsets of $X$ into bounded subsets of $\mathbb C$. The set
\[
\H_b(X)=\{f:X\to\mathbb C,\  f \textrm{ is a holomorphic function of bounded type}\}
\] is a Fr\'{e}chet algebra if we endow it with the family of seminorms $\big\{\|\cdot\|_{rB_X}:\ r>0\big\}$, where
$$
\|f\|_{rB_X}=\sup_{\|x\|<r}|f(x)|.
$$

Based on ideas due to Arens \cite{arens1951adjoint, arens1951operations}, there is a canonical way to extend an element $A\in \mathcal{L}(^mX,E)$ to an element $\widetilde{A}\in \mathcal{L}(^mX'',E'')$ using weak star limits. Given $x_1'',\ldots, x_n''\in X''$, we define
$$\widetilde{A}(x_1'',\ldots, x_n'') =  \wslim_{\,\,\,\,\,\,\,\,\,\,\,\,\,\,\,\,\,\,\,\,\,\,\alpha_1}\cdots \wslim_{\,\,\,\,\,\,\,\,\,\,\,\,\,\,\,\,\,\,\,\,\,\,\alpha_m}  A(x_{\alpha_1},\ldots, x_{\alpha_m}),$$
where $(x_{\alpha_i})\subseteq X$ is any net $w^*$-convergent to $x''_i$. In other words, we extend one variable at a time, from right to left, using  weak star limits. From now on we write
$$\wslim_{\alpha_1,\ldots, \alpha_m} A(x_{\alpha_1},\ldots, x_{\alpha_m})$$
for short. This is not the original definition (see \cite{AroBer78}), but it is equivalent to it and it is the one used through this work. Naturally, this way of extending multilinear mappings gives a tool to extend a polynomial $P \in \mathcal{P}(^mX,E)$ to  $\widetilde{P}\in \mathcal{P}(^mX'',E'')$ by defining
$$\widetilde{P}(x'')= \widetilde{\widecheck{P}}(x'',\cdots, x'').$$
The above provides a canonical extension $[f\mapsto \widetilde f]$  from $\H_b(X)$ to $\H_b(X'')$ which  satisfies, for any $f\in\H_b(X)$ and  $r>0$,
\[
\|f\|_{rB_X} = \|\widetilde f\|_{rB_{X''}}.
\] 
This extension of holomorphic functions is often referred to as the Aron-Berner extension.  We refer the reader to the articles \cite{AroBer78,DavGam89}, or the survey \cite{zalduendo2005extending} (and the references therein), for more details on this topic.

A fundamental property to consider when studying canonical extensions to the bidual is the \emph{Arens regularity}: a Banach space $X$ is said to be {\em Arens regular}, or simply \emph{regular}, if  every continuous  linear mapping $T\in\mathcal{L}(X, X')$ turns out to be weakly compact. When this holds just  for every {\em symmetric} mapping $T$
  (i. e., those satisfying $T(x_1)(x_2)=T(x_2)(x_1)$ for all $x_1, x_2\in X$) the space $X$ is said to be {\em symmetrically regular}. Note that regularity is more restrictive than symmetric regularity: it is commented in \cite{AroGalGarMae96} that the dual $J'$ of  James' space is symmetrically regular although it is not regular.



\section{Evaluations in  higher duals}\label{Sect:Cuarto dual}

As mentioned in the introduction, in \cite[Theorem 1.3]{AroGalGarMae96} it is proved that a space is non-symmetrically regular if and only if there is a fiber $\mathcal{M}_{x''}$ with an evaluation homomorphism $\tilde{\tilde{\delta}}_\cd{x}\in \mathcal{M}_{x''}$, different from $\tilde{\delta}_{x''}$. In the main result of this section we  prove something stronger for a large class of  non-symmetrically regular Banach spaces: on {\bf every} fiber there is such a homomorphism. Also, we show that for an arbitrary non-symmetrically regular space a similar result holds, eventually seeking for evaluations at points of the sixth dual. The main result of this section is the next theorem, whose proof will be  presented below.

\begin{thm}\label{teo evaluaciones cuarto}
Suppose the Banach space $X$ has  a quotient isomorphic to  $F\oplus F'$ for some non-reflexive Banach space $F$. 
Then, for every $x''\in X''$  there is $\cd{x}\in \cd{X}$ such that $$\tilde{\tilde{\delta}}_\cd{x}\in \mathcal{M}_{x''}\,\,\,\text{and}\,\,\,\tilde{\tilde{\delta}}_{\cd{x}}\neq \tilde{\delta}_{x''}.$$
\end{thm}

A large family  covered by the previous theorem is that of  Banach spaces containing complemented (isomorphic) copies of $\ell_1$. 
Indeed, such a Banach space $X$ has a quotient isomorphic to $\ell_1$. Since $c_0\oplus \ell_1$ is separable, there is a quotient map from $\ell_1$ onto $c_0\oplus \ell_1$  and therefore $X$  has a quotient isomorphic to $c_0\oplus \ell_1$. Hence, $X$ satisfies the hypotheses of the above theorem for $F=c_0$.  This gives the following.

\begin{cor}\label{coro ele 1}
Let $X$ be a Banach space  containing a complemented copy of $\ell_1$. 
Then, for every $x''\in X''$  there is $\cd{x}\in \cd{X}$ such that $$\tilde{\tilde{\delta}}_\cd{x}\in \mathcal{M}_{x''}\,\,\,\text{and}\,\,\,\tilde{\tilde{\delta}}_{\cd{x}}\neq \tilde{\delta}_{x''}.$$
\end{cor}

One may note that the argument of the previous corollary holds for any space $X$ having a quotient isomorphic to $\ell_1$, not only for those containing complemented copies of $\ell_1$. However, this observation will not give us more examples, since the existence of such a quotient is equivalent to $X$ having a complemented copy of $\ell_1$ (see for example \cite[Chapter 7, Theorem 5 ]{diestel2012sequences}).
We also mention that the hypotheses of Theorem \ref{teo evaluaciones cuarto} are more general than those of the corollary: 
the space $X= J\oplus J'$, where $J$ is James' space, satisfies the hypotheses of the theorem  but does not contain copies of $\ell_1$ since it has separable dual.

The previous results cover most of the known examples of non-symmetrically regular spaces in the literature. Obvious examples are infinite dimensional $L_1(\mu)$ spaces. 
Moreover,  if $X'$ is complemented in a Banach lattice or if it is isomorphic to a closed subspace of an ordered continuous Banach lattice, then the condition of having a complemented copy of $\ell_1$ is actually equivalent to being non-symmetrically regular (see \cite[Corollary 10]{saab1993extensions}).
In fact, we are not aware of any non-symmetrically regular Banach space not covered by the theorem. But for those possible cases, we can see that every fiber has evaluations in the sixth dual which are different from all  the evaluations in the bidual.

\begin{prop}\label{prop sexto dual}
Let $X$ be a non-symmetrically regular Banach space and $x''\in rB_{X''}$. Then, there is $\sd{x}\in rB_\sd{X}$ such that
$$\tilde{\tilde{\tilde{\delta}}}_\sd{x}\in \mathcal{M}_{x''} \,\,\,\text{ and }\,\,\,\tilde{\tilde{\tilde{\delta}}}_\sd{x}\neq \tilde{\delta}_{x''}.$$
\end{prop}
We remark that the homomorphism $\tilde{\tilde{\tilde{\delta}}}_\sd{x}$ given by the previous proposition may actually be an evaluation in $X^{iv}$, in the sense that there  may be a $\cd{y}\in \cd{X}$ such that  $\tilde{\tilde{\tilde{\delta}}}_\sd{x}=\tilde{\tilde{\delta}}_\cd{y}.$
\medskip

Now we focus on the proof of Theorem \ref{teo evaluaciones cuarto}. We begin by proving the following lemma, which gives sufficient conditions for the existence of  evaluations in the fourth dual giving new homomorphisms in every fiber of the spectrum.

\begin{lem}
\label{prop 1} If there exist an element $w\in (X')^\bot\subseteq X^{iv}$ and a polynomial $P\in \mathcal{P}(^mX)$ satisfying
$$\widetilde{\widetilde{P}}(w)\neq 0$$
then, for any $x''\in X''$ there is $\cd{x}\in X^{iv}$ such that
$$\tilde{\tilde{\delta}}_{x^{iv}}\in \mathcal{M}_{x''}\,\,\,\text{and}\,\,\,\tilde{\tilde{\delta}}_{\cd{x}}\neq \tilde{\delta}_{x''}.$$
\end{lem}
\begin{proof} Recall that $X^{iv}=X''\oplus (X')^\bot$ and
take $\cd{x}=x''+\alpha w$ with $\alpha \in \mathbb{C}$ to be fixed. It is clear that $\tilde{\tilde{\delta}}_{\cd{x}}\in \mathcal{M}_{x''}$ and we can compute
$$\widetilde{\widetilde{P}}(\cd{x})=\widetilde{\widetilde{P}}(x''+\alpha w)=\sum_{k=0}^m {{m}\choose{k}} \widecheck{\widetilde{\widetilde{P}}}((x'')^{m-k},(\alpha w)^k)=\sum_{k=0}^m {{m}\choose{k}} \alpha^k\widecheck{\widetilde{\widetilde{P}}}((x'')^{m-k},w^k),$$
obtaining a polynomial of degree $m$ on $\alpha$. Thus, we can set $\alpha$ such that 
$$\tilde{\tilde{\delta}}_{\cd{x}}(P)=\widetilde{\widetilde{P}}(x''+\alpha w) \neq \widetilde{\widetilde{P}}(x''+0 w) =\widetilde{P}(x'')=\tilde{\delta}_{x''}(P). \qedhere$$ 
\end{proof}

Note that the hypothesis of the previous lemma is equivalent  to the existence of elements $x_0\in X$, $w\in (X')^\bot$ and a polynomial (non necessarily homogeneous) $Q$ such that $\widetilde{\widetilde{Q}}(x_0+ w)\neq Q(x_0)$.

\bigskip

Let $P\in \mathcal{P}(^2X)$ be a polynomial and $A$ be its associated symmetric bilinear form. We define the \emph{linear operator associated to $P$ } as  $T:X\to X'$  given by $$T(x)= A(x,\cdot).$$
Note that the linear operator $T$ associated to a polynomial $P$ is symmmetric. Moreover, by \cite[Lemma 8.1]{AroColGam91} the extension of $P$ can be computed as
\begin{equation}\label{eq-extension-transpuesta} 
\widetilde{P}(x'')= \langle T'' x'', x''\rangle.
\end{equation}

\begin{rem}\label{rem-ext4to}
Let $P\in \mathcal{P}(^2X)$ with associated operator $T\in\mathcal L(X,X')$. 
By \cite[Lemma 8.1]{AroColGam91},  the symmetric operator associated to 
$\widetilde{P}$ is $T_{\widetilde{P}}\in\mathcal L(X'', X''')$ given by $T_{\widetilde{P}}=\frac{T''+i_{X'}i_X' T''}{2}$. Now we can use \eqref{eq-extension-transpuesta} to compute the extension $\widetilde{\widetilde{P}}$ of $P$ to $\cd{X}$ as 
\begin{eqnarray*}
\widetilde{\widetilde{P}}(\cd{x})&=&\left\langle T_{\widetilde{P}}''\cd{x},\cd{x}\right\rangle\\
&=& \left\langle \left(\frac{T''+i_{X'}i_X' T''}{2}\right)''(\cd{x}),\cd{x}\right\rangle\\
&=&\left\langle  \frac{\cd{T}(\cd{x})+i_{X'}''i_X''' \cd{T}(\cd{x})}{2},\cd{x}\right\rangle.\
\end{eqnarray*}
\end{rem}

During the course of the present work, we noticed that in the proof of \cite[Theorem 1.3]{AroGalGarMae96} the extension $\widetilde{\widetilde{P}}$ is not properly computed and subsequent steps are then not valid. Fortunately, using the expression for $\widetilde{\widetilde{P}}$ from  the previous remark,  that  proof  can be fixed with some minor changes that do not alter the original idea. Since that result is cited repeatedly here, for the sake of completeness and the convenience of the reader, we specify these changes below. We are grateful to Richard Aron and Manolo Maestre for useful talks regarding this topic.

In \cite[Theorem 1.3]{AroGalGarMae96}, the authors need to show that the linear coefficient of the polynomial of degree 2
\begin{equation}\label{eq-that-polynomial}
    \alpha\longrightarrow \widetilde{\widetilde{P}}(x_0''+\alpha w)
\end{equation}
is not null. Here, $P$ is a $2-$homogeneous polynomial in $X$ whose associated operator $T$ is not weakly compact, and $x_0'' \in X'', w\in (X')^\bot \subseteq X^{iv}$ are chosen such that \linebreak $\langle T''x_0'' , w \rangle \neq 0$. The existence of these elements is a consequence of $T$ not being weakly compact. The linear coefficient of the polynomial given in \eqref{eq-that-polynomial} is
\begin{equation}\label{linear coefficient} \frac{1}{2}\left(\langle T^{iv}x_0'' , w \rangle + \langle i_{X'}''i_X''' T^{iv} x_0'', w \rangle  + \langle T^{iv}w, x_0''\rangle + \langle i_{X'}''i_X'''T^{iv} w ,x_0'' \rangle \right).\end{equation}
The term $\langle T^{iv}x_0'' , w \rangle$ is equal to $\langle T''x_0'' , w \rangle$.
For the second term we have that
$$ \langle i_{X'}''i_X''' T^{iv} x_0'', w \rangle =  \langle i_{X'}i_X' T'' x_0'', w \rangle=  \langle i_X' T'' x_0'',i_{X'}' w \rangle$$
which is null since $i_{X'}' w=0$ (recall that $w\in (X')^\bot$).
Similarly 
$$\langle T^{iv}w, x_0''\rangle  = \langle w,T'''x_0''\rangle =\langle w,T'x_0''\rangle =0$$
 because $T'x_0''\in X'$.
Finally, using that $T$ is symmetric we have that $T''=T'''i_X''$, and thus
$$\langle i_{X'}''i_X'''T^{iv} w ,x_0'' \rangle =\langle  w , T'''i_X''i_{X'}'x_0'' \rangle= \langle  w , T''x_0'' \rangle$$
which is exactly the first term of \eqref{linear coefficient}.
Therefore, the linear coefficient \eqref{linear coefficient} is equal to $\langle T''x_0'' , w \rangle$, which incidentally is the same obtained by the authors in \cite{AroGalGarMae96}, despite the computation mistake.

\bigskip
The following example is crucial for the proof of the main theorem of this section. It shows that in any space of the form $F\oplus F'$, with $F$ non-reflexive, we can find in every fiber evaluations at elements of the fourth dual giving new homomorphisms.

 \begin{exa}\label{exa evaluaciones} Given a non-reflexive Banach space $F$, consider $E= F\oplus F'$  and $P\in \mathcal{P}(^2E)$ defined as 
 $$P((x,x'))=x'(x).$$ Then, $P$ satisfies the conditions of Lemma \ref{prop 1}. As a consequence, in every fiber there are evaluations at elements of the fourth dual giving new homomorphisms.
\end{exa}
To see this, we first note that $\widecheck{P}((x,x'),(y,y'))=\frac{x'(y)+y'(x)}{2}$, so the operator \linebreak $T:F\oplus F' \rightarrow F'\oplus F''$ associated to $P$ is given by
 $$T((x,x'))=\frac{1}{2}(x', i_F(x)).$$
 Therefore, by Remark \ref{rem-ext4to} we have that $\widetilde{\widetilde P}\in\mathcal P(^2 E^{iv})$ is computed in the following way:
 \begin{eqnarray*}
\widetilde{\widetilde{P}}((\cd{x},x^v))&=& \left\langle  \frac{\cd{T}((\cd{x},x^v))+i_{E'}''i_E''' \cd{T}((\cd{x},x^v))}{2},(\cd{x},x^v)\right\rangle\\
&=&\left\langle  \frac{(x^v,i_F^{iv}\cd{x}) +  i_{E'}''i_E'''(x^v,i_F^{iv}\cd{x}))}{4},(\cd{x},x^v)\right\rangle \\
&=&\left\langle  \frac{(x^v,i_F^{iv}\cd{x}) +  (i_{F'}''i_F'''(x^v), i_{F''}''i_{F'}'''i_F^{iv}\cd{x})}{4},(\cd{x},x^v)\right\rangle \\
&=& \frac{x^v(\cd{x})+i_F^{iv}\cd{x}(x^v)+i_{F'}''i_F'''(x^v)(\cd{x}) + i_{F''}'' id_{F'}''' \cd{x}(x^v)}{4}\\
&=& \frac{x^v(\cd{x})+i_F^{iv}\cd{x}(x^v)+i_F'''(x^v)(i_{F'}'\cd{x}) + \cd{x}(i_{F''}'x^v)}{4}.\
 \end{eqnarray*}

 Now, if we consider that $(\cd{x},\qd{x})$ belongs to $(F')^\bot\oplus (F'')^\bot=  (E')^\bot\subseteq \cd{E}$ with $\qd{x}\not=0$ we have $i_{F'}'\cd{x}=0$ and $i_{F''}'\qd{x}=0$. Thus,
 $$\widetilde{\widetilde{P}}((\cd{x},\qd{x}))=\frac{\qd{x}(\cd{x})+i_F^{iv}\cd{x}(\qd{x})}{4}=\frac{\qd{x}(\cd{x})+\cd{x}(i_F'''\qd{x})}{4}.$$
 Recall that $F'''=F'\oplus F^\bot$ and let us call $\rho:F'''\rightarrow F^\bot$ the canonical projection.  The fact that $\cd{x}$ belongs to $(F')^\bot$ implies $\cd{x}=\cd{x}\rho$. Therefore, 
 \begin{eqnarray*}
 \widetilde{\widetilde{P}}((\cd{x},\qd{x}))&=&\frac{\qd{x}(\cd{x})+\cd{x}(\rho i_F'''\qd{x})}{4} = \frac{\qd{x}(\cd{x})+i_{F'''}\rho i_F'''\qd{x} (\cd{x})}{4}\\
 &=& \frac{[\qd{x}+i_{F'''}\rho i_F'''\qd{x}] (\cd{x})}{4}.\
 \end{eqnarray*}

Note that $\qd{x}+i_{F'''}\rho i_F'''\qd{x}$ cannot belong to the copy of $F'$ inside $F^{v}$ due to the next observations:
 \begin{itemize}
     \item $F^v=F'''\oplus (F'')^\bot=F'\oplus F^\bot\oplus (F'')^\bot$
     \item $\qd{x}\in (F'')^\bot\setminus \{0\}$
     \item $i_{F'''}\rho i_F'''\qd{x}\in F^\bot.$
 \end{itemize}
Then, we can choose $\cd{x}\in (F')^\bot$ such that $[\qd{x}+i_{F'''}\rho i_F'''\qd{x}] (\cd{x})\neq 0$, reaching in this way to an element $(\cd{x},\qd{x})\in (E')^\bot$ that satisfies $\widetilde{\widetilde{P}}((\cd{x},\qd{x}))\neq 0$. Finally, appealing to Lemma \ref{prop 1} we get the desired conclusion.
 
\bigskip
In \cite{leung1993banach} there is an example of a non-reflexive Banach $F$ space such that both $F$ and $F'$ are regular (and in particular simetrically regular). Even then, the space $F\oplus F'$ not only is non-symmetrically regular (see also the comments below \cite[Corollary 1.8]{AroGalGarMae96}), but moreover in every fiber there are evaluations in the fourth dual giving new homomorphisms.

\bigskip 

We are ready to prove the main theorem of this section.

\begin{proof}[Proof of Theorem \ref{teo evaluaciones cuarto}] Writing $E=F\oplus F'$,
Example \ref{exa evaluaciones} gives us a polynomial $Q:E\rightarrow \zC$ with $\widetilde{\widetilde{Q}}(u)\neq 0$ for some $u\in  (E')^\bot$.  Now, the assumptions on $X$ mean that there is a continuous surjective linear map $T:X\to E$.
Let us see that $P:=Q\circ T$ satisfies the conditions of Lemma \ref{prop 1}. Using that
$$\widetilde{\widetilde{P}}=\widetilde{\widetilde{Q}}\circ \cd{T}$$
and that $\cd{T}$ is onto, we can take $\cd{x}\in X^{iv}$ such that $\cd{T}(\cd{x})=u$. In particular
$$\widetilde{\widetilde{P}}(\cd{x})=\widetilde{\widetilde{Q}}(u)\neq 0.$$
Let us write $\cd{x}=x''+w$, with $x''\in X''$ and $w\in (X')^\bot$. Then
$$u=\cd{T}(\cd{x})=T''(x'')+\cd{T}(w).$$
Since
\begin{itemize}
    \item $w\in (X')^\bot$ implies that $\cd{T}(w)\in (E')^\bot$
    \item $T''(x'')\in E''$
    \item $u\in (E')^\bot$
\end{itemize}
we conclude that $T''(x'')=0.$ Therefore, $\widetilde{\widetilde{P}}(w)= \widetilde{\widetilde{P}}(\cd{x})\ne 0$ and given that $w\in (X')^\bot$ we can use Lemma \ref{prop 1} to finish the proof.
\end{proof}

We continue with the proof of Proposition \ref{prop sexto dual}.

\begin{proof}[Proof of Proposition \ref{prop sexto dual}]
Let $P$ be a $2$-homogeneous polynomial such that its associated operator $T$ is not weakly compact. Since $T$ factors through   $T_{\widetilde{P}}$ (the operator associated to $\widetilde{P}$), we have that  $T_{\widetilde{P}}$ cannot be weakly compact. Thus, there is $\cd{x}\in \cd{X}$ such that
$$T_{\widetilde{P}}''(\cd{x})\in X^v\setminus X'''.$$

Since 
$$\cd{X}= X''\oplus X^\bot\,\,\,\text{ and }\,\,\, T_{\widetilde{P}}''(X'')=T_{\widetilde{P}}(X'')\subseteq X''',$$
eventually replacing $\cd{x}$ by $\cd{x}-i_{X'}'(\cd{x})+x''$, we may assume $i_{X'}'(\cd{x})=x''$. Replacing $\cd{x}$ by $\frac{\cd{x}-x''}{K}+x''$, with $K>0$ large enough, we may also assume that $\Vert \cd{x} \Vert<r$.

Now, the result follows using a similar argument as in the proof  of \cite[Theorem 1.3]{AroGalGarMae96}. Since $T_{\widetilde{P}}''(\cd{x})\in X^v\setminus X'''$, there is $w\in (X''')^\bot\subseteq \sd{X}$ such that 
$$\alpha \rightarrow \widetilde{\widetilde{\widetilde{P}}}(\cd{x}+\alpha w)$$
is a polynomial with non-null linear coefficient (see the comments after Remark \ref{rem-ext4to}). In particular, there is $\alpha_0$ close to $0$, such that
$$\Vert \cd{x}+\alpha_0 w\Vert <r\,\,\,\text{ and } \,\,\, \widetilde{\widetilde{\widetilde{P}}}(\cd{x}+\alpha_0w)\neq \widetilde{P}(x'').$$
Taking $\sd{x}=\cd{x}+\alpha_0w$ we conclude the proof.
\end{proof}

\subsection*{A remark on the corona of the spectrum}

Throughout this article we show some similarities and differences between the evaluations at points in the bidual and evaluations at points in  higher duals. We recall that, in our setting, the corona is defined as the set $\M_b(X)\setminus \overline{\{\delta_x: x\in X \}} $ (the closure taken in the Gelfand topology), i.e., the set of homomorphisms that cannot be approximated by evaluations in $X$.
We end this section with the observation that the results in \cite{DavGam89}, which clearly imply that  evaluations at points of  the bidual do not belong to the corona of the spectrum, also imply that evaluations at points of higher duals do not belong either. 

Note that most of the examples of homomorphisms appearing in the literature  are given by evaluations at points of some  dual $X^{(2n)}$ or by limits of evaluations along ultrafilters (see for example \cite[Proposition 1.5]{AroGalGarMae96}). The latter are clearly outside the corona, so it is interesting to determine what happens to evaluations in $X^{(2n)}$.

\begin{prop}\label{prop-corona} Let $X$ be an arbitrary Banach space, then
$$\{\tilde{\tilde{\delta}}_\cd{x},\  \cd{x}\in \cd{X} \} \subseteq \overline{\{\delta_x,\  x\in X \}},$$
where the closure is taken in the Gelfand topology.
\end{prop}

\begin{proof}
Consider the sets
\begin{align*} 
\Delta &= \{\delta_x,\ x\in X \}\subseteq \mathcal M_b(X) \\
\Delta_2 &= \{\tilde{\delta}_{x''},\  x''\in X'' \}\subseteq \mathcal M_b(X)\\
\Delta_4 &= \{\tilde{\tilde{\delta}}_\cd{x},\  \cd{x}\in \cd{X} \}\subseteq \mathcal M_b(X)\\
\Theta_2 &= \{\delta_{x''},\  x''\in X \}\subseteq \mathcal M_b(X'')\\
\Theta_4 &= \{\tilde{\delta}_\cd{x},\  \cd{x}\in \cd{X} \}\subseteq \mathcal M_b(X'').\
\end{align*}
By \cite{DavGam89} we know that $\Delta_2\subseteq \overline{\Delta}$ and $\Theta_4\subseteq \overline{\Theta_2}$. The result follows by noticing that $\Theta_2$ and $\Theta_4$ naturally map, respectively, onto $\Delta_2$ and $\Delta_4$ in $\mathcal M_b(X)$. 

More precisely, consider the function $\varrho: \mathcal M_b(X'')\rightarrow \mathcal M_b(X)$ defined as $\varrho(\varphi)(f)=\varphi(\tilde{f})$. Since $\varrho$ in continuous, we have that  $\varrho(\overline{\Theta_2})\subseteq \overline{\varrho(\Theta_2)},$ thus
$$ \Delta_4 = \varrho(\Theta_4)\subseteq \varrho(\overline{\Theta_2})\subseteq \overline{\varrho(\Theta_2)}=\overline{\Delta_2}\subseteq \overline{\overline{\Delta}}=\overline{\Delta}. \qedhere$$
\end{proof}

It is not hard to adapt the proof above to show the same result for evaluations in higher even duals.



\section{The convolution operation}\label{Sect:convolucion}

In this section we study the convolution operation on $\mathcal M_{b}(X)$. We extend or complement some results on commutativity from \cite[Section 6]{AroColGam91}, with a special focus on the non-symmetrically regular case. 
Recall that
given $\varphi $, $\psi \in \mathcal M_{b}(X)$ we can define their
convolution $ \varphi \ast \psi $ as
\begin{equation*}
\varphi \ast \psi (f)=\varphi (x\mapsto \psi (\tau
_{x}(f)))\text,
\end{equation*}
where $\tau
_{x}$ is the translation operator: $\tau
_{x}(f)(y)=f(y+x)$. 

Convolution is not commutative in general. However, we know from \cite[Lemma 6.3]{AroColGam91} that $\delta_x\ast\varphi=\varphi\ast\delta_x$ for all $x\in X$ and all $\varphi \in \mathcal M_{b}(X)$. Regarding evaluations in the bidual, 
we state as a theorem a combination of Lemma 6.10 and Theorem 6.11 from \cite{AroColGam91}.

\begin{thm}\cite{AroColGam91}
\label{lem conv que conmutan} Given $z_1, z_2\in X''$, the following statements are equivalent. 
\begin{enumerate}
    \item[(i)] $\tilde\delta_{z_1} \ast \tilde\delta_{z_2} =\tilde\delta_{z_2} \ast \tilde\delta_{z_1}$.
    \item[(ii)] $\tilde\delta_{z_1} \ast \tilde\delta_{z_2} =\tilde\delta_{z_1+z_2}$.   
    \item[(iii)]  For any $P\in \mathcal{P}(^2X)$ we have $
\widetilde{\widecheck{P}}(z_1,z_2)= \widetilde{\widecheck{P}}(z_2,z_1).$
\end{enumerate}
\end{thm}
As a consequence of this result, all evaluations in the bidual commute if and only if the space $X$ is symmetrically regular. This means that in every non-symmetrically regular Banach space $X$ there exists a pair $z_1,z_2\in X''$ such that $\tilde\delta_{z_1}$ and $\tilde\delta_{z_2}$ do not commute. The main goal of this section is to show that $X=\ell_1$ and $X=F\oplus F'$ for a non-reflexive $F$ are extreme cases in the following sense: the only evaluations that commute with every other evaluation are those at points in $X$.

\begin{thm}\label{teo main sec2}
Let $X$ be $\ell_1$ or $F\oplus F'$, with $F$ a non-reflexive Banach space and take  $z\in X''$. If  $\tilde\delta_z\ast\tilde\delta_u=\tilde\delta_u\ast \tilde\delta_z$ for all $u\in X''$, then necessarily  $z\in X$.
\end{thm}
Combining the previous theorem with Proposition \ref{cor morfismos que conm} below, we also see that if a homomorphism $\varphi$ commutes with any other homomorphism (or any other evaluation in the bidual), then $\pi(\varphi)$ must belong to $X$.

In order to prove Theorem \ref{teo main sec2}, we first deal with $X=\ell_1$. The fact that the non-symmetrically regular space $\ell_1$ is itself a dual space yields the decomposition \linebreak $\ell_1''=\ell_1\oplus c_0^\bot$ which is relevant to our arguments.
Next lemma provides us of a family of symmetric bilinear forms on $\ell_1$ which behave as product of evaluations  at some elements in $\ell_1''$.

\begin{lem}
\label{prop bilineal de ele1}
Given $a,b\in \ell_\infty$, there is a symmetric bilinear form
$B\in\mathcal L(^2\ell_1)$
such that, if $z_1\in\ell_1''$ and $z_2\in c_0^\bot\subseteq \ell_1''$, then
$$\widetilde{B}(z_1,z_2)=z_1(a)z_2(b).$$

\end{lem}
\begin{proof}
If $a=(a_n)_{n\in \zN}$ and $b=(b_n)_{n\in\zN}$, given $k\in \zN$ let us call 
$$a^{\leq k}=(a_1,\ldots, a_k,0,\ldots).$$
$$b^{\geq k}=(0,\ldots, 0,b_{k},b_{k+1}, \ldots).$$
Define the bilinear form $\phi\in\mathcal L(^2\ell_1)$ as 
\begin{eqnarray*}
\phi(x,y)&=&\sum_{n= 1}^\infty \sum_{k=1}^n b_nx_na_ky_k\\
&=&\sum_{n= 1}^\infty b_nx_na^{\leq n}(y)\\
&=&\sum_{k= 1}^\infty a_ky_kb^{\geq k}(x).\
\end{eqnarray*}
Let us show that the symmetric bilinear form 
$$B(x,y)=\phi(x,y)+\phi(y,x)$$
meets the require conditions. Let $z_1\in\ell_1''$ and $z_2\in c_0^\bot\subseteq \ell_1''$ and observe that $a^{\leq n}, (b-b^{\geq k})\in c_0$, which implies that
$$z_2(a^{\leq n})=0  \,\,\,\text{ and }  z_2(b^{\geq k})=z_2(b).$$
Then, if $(x^\alpha),(y^\beta)\subseteq \ell_1$ are nets $w^*$-convergent to $z_1$ and $z_2$ respectively, we have 
\begin{eqnarray*}
\widetilde{B}(z_1,z_2)&=& \lim_{\alpha, \beta} B(x^\alpha, y^\beta) = \lim_{\alpha, \beta} \phi(x^\alpha, y^\beta)+\phi(y^\beta, x^\alpha)\\
&=& \lim_{\alpha, \beta} \left(\sum_{n=1}^\infty b_n \,x^\alpha_n \,\,a^{\leq n}(y^\beta) + \sum_{k=1}^\infty a_k \,x^\alpha_k \,\, b^{\geq k}(y^\beta)\right)\\
&=& \lim_{\alpha} \left( \sum_{n=1}^\infty b_n \,x^\alpha_n \,\, z_2(a^{\leq n}) + \sum_{k=1}^\infty a_k \,x^\alpha_k  \,\, z_2(b^{\geq k}) \right) \\
&=& \lim_{\alpha}  \left( 0 + \sum_{k=1}^\infty a_k x^\alpha_k z_2(b) \right)= \lim_{\alpha} z_2(b) a(x^\alpha)\\
&=& z_2(b)z_1(a).\qedhere
\end{eqnarray*}
\end{proof}

The previous lemma allows us to give a more direct proof of \cite[Theorem 7.5]{AroColGam91} without appealing to measure theoretic tools developed in that article. 

\begin{thm}\label{thm ACG}\cite[Theorem 7.5]{AroColGam91}
Given nonzero $z_1,z_2\in c_0^\bot\subseteq \ell_1''$, then $\tilde{\delta}_{z_1}\ast \tilde{\delta}_{z_2}=\tilde{\delta}_{z_2}\ast \tilde{\delta}_{z_1}$ if and only if $z_1$ and $z_2$ are multiples of each other. 
\end{thm}
\begin{proof}
If $z_1$ and $z_2$ are not multiples of each other, there are $a,b\in \ell_\infty$ such that 
$$a\in \ker(z_2)\setminus \ker(z_1)\,\,\,\text{ and }\,\,\, b\notin\ker(z_2).$$
For this $a$ and $b$, take $B$ as in Lemma \ref{prop bilineal de ele1}. Thus
$$\widetilde{B}(z_1,z_2)=z_1(a)z_2(b)\neq 0 =z_2(a)z_1(b)=\widetilde{B}(z_2,z_1).$$
Then, the result follows by Theorem \ref{lem conv que conmutan}. 
\end{proof}

Next corollary is an extension of the previous theorem which will be useful for the proof of our main theorem. 

\begin{cor}\label{cor comn de ele1}
Given $z_1,z_2\in \ell_1''$. If
$z_1=x_1+w_1$ and $z_2=x_2+w_2,$
with $x_1,x_2\in \ell_1$ and $w_1,w_2\in c_0^\bot\setminus\{0\} $, then  $\tilde{\delta}_{z_1}\ast \tilde{\delta}_{z_2}=\tilde{\delta}_{z_2}\ast \tilde{\delta}_{z_1}$ if and only if $w_1$ and $w_2$ are multiples of each other.
\end{cor}
\begin{proof}
Since evaluations at points of $\ell_1$ always commute, by Theorem \ref{lem conv que conmutan}, we have  
\begin{eqnarray*}
\tilde{\delta}_{z_1}\ast \tilde{\delta}_{z_2}&=& \tilde{\delta}_{x_1+w_1}\ast \tilde{\delta}_{x_2+w_2}\\
&=& \delta_{x_1}\ast\tilde{\delta}_{w_1}\ast \delta_{x_2}\ast\tilde{\delta}_{w_2}\\
&=&\delta_{x_1}\ast\delta_{x_2}\ast\tilde{\delta}_{w_1}\ast\tilde{\delta}_{w_2}.\\
&=&\delta_{x_1+x_2}\ast\tilde{\delta}_{w_1}\ast\tilde{\delta}_{w_2}.\
\end{eqnarray*}
Similarly
\begin{equation*}
    \tilde{\delta}_{z_2}\ast \tilde{\delta}_{z_1}=\delta_{x_2+x_1}\ast\tilde{\delta}_{w_2}\ast\tilde{\delta}_{w_1}.\
\end{equation*}
Then, since  $\delta_{x_1+x_2}$ is invertible, $\tilde{\delta}_{z_1}\ast \tilde{\delta}_{z_2}=\tilde{\delta}_{z_2}\ast \tilde{\delta}_{z_1}$ if and only if    $\tilde{\delta}_{w_1}\ast \tilde{\delta}_{w_2}=\tilde{\delta}_{w_2}\ast \tilde{\delta}_{w_1}$. Finally, an appeal to Theorem \ref{thm ACG} finishes the proof.
\end{proof}

We now focus on the case $X=F\oplus F'$. The main step for the proof is the next lemma which is a kind of  analogue of  Corollary \ref{cor comn de ele1}.

\begin{lem}\label{lem evaluaciones}
Fix $x'',y''\in F'' $ and $x''',y'''\in F'''$. If $\tilde{\delta}_{(x'',x''')}$ and $\tilde{\delta}_{(y'',y''')}$ commute,  then 
 $\rho(x''')(y'') = \rho(y''')(x'')$, where $\rho:F''' \rightarrow F^\bot$ is the canonical projection. 
\end{lem}

\begin{proof} 
Consider  the polynomial $P$ from Example \ref{exa evaluaciones}. By Theorem \ref{lem conv que conmutan} $$\widetilde{\widecheck{P}}((x'',x'''),(y'',y'''))=\widetilde{\widecheck{P}}((y'',y'''),(x'',x''')).$$ Writing $x'''=i_{X}'x'''+\rho (x''')$ we obtain
\begin{eqnarray*}
\widetilde{\widecheck{P}}((x'',x'''),(y'',y''')) &=& \frac{x'''(y'')+x''(i_X' y''')}{2}\\
&=& \frac{y''(i_{X}'x''')+\rho(x''') (y'')+x''(i_X' y''')}{2}.\
\end{eqnarray*}
Similarly 
$$\widetilde{\widecheck{P}}((y'',y'''),(x'',x'''))=\frac{x''(i_X' y''')+\rho(y''') (x'')+y''(i_{X}'x''')}{2},$$
which gives the desired conclusion.
\end{proof}

\begin{proof}[Proof of Theorem \ref{teo main sec2}] We start with the case $X=\ell_1$.
We write $z=x+w$ with $x\in \ell_1$ and $w\in c_0^\bot$. 
If $w\ne 0$, by Corollary~\ref{cor comn de ele1}, $\tilde\delta_z$ commutes with  $\tilde\delta_u$ only if the projection of $u$ on $c_0^\bot$ is a multiple of $w$. Therefore, if $\tilde\delta_z$  commutes with every other evaluation in the bidual we must have $w=0$ or, equivalently, $z\in \ell_1$.

The argument for the case $X= F\oplus F'$ is similar  using Lemma~\ref{lem evaluaciones} instead of Corollary~\ref{cor comn de ele1}.
\end{proof}

We finish this section by stating and proving the proposition mentioned after Theorem~\ref{teo main sec2}.

 \begin{prop}\label{cor morfismos que conm} Let $X$ be  a Banach space. Given two homomorphisms $\varphi\in \mathcal{M}_{z_1}$ and  $\psi\in \mathcal{M}_{z_2}$, if $\varphi \ast \psi = \psi \ast \varphi$ then $\tilde{\delta}_{z_1} \ast \tilde{\delta}_{z_2}= \tilde{\delta}_{z_2}\ast \tilde{\delta}_{z_1}$. 
 \end{prop}

 \begin{proof}
Let $P$ be a $2-$homogeneous polynomial. Then
\begin{equation}\label{ec conm 1}
\psi\big(\tau_{x}(P)\big) = \psi\big(P(\,\cdot \,)+2\widecheck{P}(x,\, \cdot \,) + P(x)\big)  = \psi(P)+ 2\widetilde{\widecheck{P}}(x,z_2)+ P(x).
\end{equation}
Therefore, $$\varphi \ast \psi(P)=\varphi(\psi(P)+ 2\widetilde{\widecheck{P}}(\,\cdot\, , z_2)+ P(\,\cdot\,))= \psi(P)+2\widetilde{\widecheck{P}}(z_1, z_2)+\varphi(P).$$
Similarly 
\begin{equation}\label{ec conm 2}
    \psi \ast \varphi(P)=\varphi(P)+2\widecheck{P}(z_2, z_1)+\psi(P).
\end{equation}
By hypothesis $\eqref{ec conm 1}$ and $\eqref{ec conm 2}$ are equal, hence $\widecheck{P}(z_1, z_2)=\widecheck{P}(z_2, z_1)$. Since this holds for any  $2$-homogeneous polynomial, by Theorem \ref{lem conv que conmutan} we conclude that $$\tilde{\delta}_{z_1} \ast \tilde{\delta}_{z_2}= \tilde{\delta}_{z_2}\ast \tilde{\delta}_{z_1}. \qedhere$$
\end{proof}

\begin{cor}
Let $X$ be $\ell_1$ or $F\oplus F'$, with $F$ a non-reflexive Banach space.  Then, the center $\mathcal{Z}$ of the semigroup $(\M_b(X), \ast)$ satisfies
$$\{ \delta_x,\  x\in X\}\subseteq \mathcal{Z} \subseteq \bigcup_{x\in X} \M_x$$
\end{cor}

Looking at Theorem \ref{lem conv que conmutan}, it is natural to ask if the statements are equivalent for evaluations in the fourth dual. Next example shows that they are not. 

\begin{exa}
Let $X$ be $\ell_1$ or $F\oplus F'$, with $F$ a non-reflexive Banach space. Then, there exists $w\in \cd{X}$ such that $\tilde{\tilde{\delta}}_{w}$ and  $\tilde{\tilde{\delta}}_{-w}$ commute but $\tilde{\tilde{\delta}}_{w}\ast \tilde{\tilde{\delta}}_{-w}\ne \delta_0$.
\end{exa}

\noindent Indeed, by the proof of Theorem \ref{teo evaluaciones cuarto} we can take $w\in \cd{X}$ and a 2-homogeneous polynomial $P$ such that $\tilde{\tilde{\delta}}_{w}\in \mathcal{M}_0$ and $\widetilde{\widetilde{P}}(w)\neq 0$. Using  \cite[Corollary 6.7]{AroColGam91} it is not hard to see that $\tilde{\tilde{\delta}}_{w}$ and  $\tilde{\tilde{\delta}}_{-w}$ commute. Let us show that $\tilde{\tilde{\delta}}_{w}\ast \tilde{\tilde{\delta}}_{-w}\ne 0$. Since $\tilde{\tilde{\delta}}_{-w}\in \mathcal{M}_0$ for each $x\in X$ we have
$$\tilde{\tilde{\delta}}_{-w}\big(\tau_{x}(P)\big) = \tilde{\tilde{\delta}}_{-w} \big(P(\,\cdot \,)+2\widecheck{P}(x,\, \cdot \,) + P(x)\big)  = \widetilde{\widetilde P}(-w) + P(x) =\widetilde{\widetilde P}(w) + P(x). $$
Therefore, $\tilde{\tilde{\delta}}_{w}\ast \tilde{\tilde{\delta}}_{-w}(P)=\tilde{\tilde{\delta}}_{w}(\widetilde{\widetilde P}(w) + P(\,\cdot\,)) =\widetilde{\widetilde P}(w)+\widetilde{\widetilde P}(w)\ne 0$.





\section{Symmetric multilinear extensions}\label{Sect:AB ext}

On symmetrically regular Banach spaces, the canonical extension of each symme\-tric bilinear mapping is again symmetric. This property is transferred to $m$-linear mappings due to the extension procedure which is on one variable at a time \cite[Section 8]{AroColGam91}. On non-symmetrically regular spaces, it is natural to search for conditions on a symmetric multilinear mapping that ensure that it has a symmetric extension. In the next result we exhibit a series of equivalences of this fact (note that some of them have a certain \textit{bilinear flavor}).

\begin{thm}\label{th AB sim} 
Given $m\ge 2$ and $A\in \mathcal{L}_s(^mX;E)$, consider its canonical extension $\widetilde{A}\in \mathcal{L}(^mX'';E'')$. The following statements are equivalent
\begin{enumerate}
    \item[(i)] $\widetilde{A}$ is symmetric.
    \item[(ii)] $\widetilde{A}$ is $(w^*-w^*)$-continuous on each variable.
    \item[(iii)]  $i_{E'}'\circ \widetilde{\widetilde{A}}=\widetilde{A}\circ (i_{X'}'\times\cdots\times i_{X'}')$.
    \item[(iv)] $\widetilde{A}$ is $(w^*-w^*)$-continuous on the second variable.
    \item[(v)] The first two variables of $\widetilde{A}$ commute.
\end{enumerate}
\end{thm}

\begin{proof}$(i)\Rightarrow (ii)$ The extension $\widetilde{A}$ is always  $(w^*-w^*)$-continuous on the first variable. Thus, if it is symmetric, it is also $(w^*-w^*)$-continuous on every variable.

$(ii)\Rightarrow (iii)$ This can be proved by adapting the argument from \cite[Lemma 1.2]{AroGalGarMae96}.  Take $z_1,\ldots, z_{m-1}\in X''$ and $x^{iv}_m\in \cd{X}$. If $(x''_\alpha) \subset X''$  is a net  converging to $x^{iv}_m$ in the $\sigma(\cd{X},X''')$-topology , then for every $e'''\in E'''$,
$$\widetilde{\widetilde{A}}(z_1,\ldots, z_{m-1},x^{iv}_m)(e''')= \lim_\alpha e'''(\widetilde{A}(z_1,\ldots, z_{m-1},x''_\alpha)).$$
This implies that, for every $e'\in E'$,
$$i_{E'}'\circ \widetilde{\widetilde{A}}(z_1,\ldots, z_{m-1},x^{iv}_m)(e')= \lim_\alpha \widetilde{A}(z_1,\ldots, z_{m-1},x''_\alpha)(e').$$

Similarly,  since $\lim_\alpha e'''(x''_\alpha)=x^{iv}_m(e''')$, for all $e'''\in E'''$ we have $\lim_\alpha x''_\alpha(e')=i_{X'}'(x^{iv}_m)(e')$ for every $e'\in E'$, meaning that $(x''_\alpha)$ is $w^*$-convergent to $i_{X'}'(x^{iv}_m)$. Now, by $(ii)$,  for every $e'\in E'$,
\begin{eqnarray*}
i_{E'}'\circ \widetilde{\widetilde{A}}(z_1,\ldots, z_{m-1},x^{iv}_m)(e')&=&\lim_\alpha \widetilde{A}(z_1,\ldots, z_{m-1},x''_\alpha)(e') \\
&=& \widetilde{A}\left(z_1,\ldots, z_{m-1},i_{X'}'(x^{iv}_m)\right)(e'). \
\end{eqnarray*}

To prove that for any $x^{iv}_1,\ldots, x^{iv}_m\in \cd{X}$
$$i_{E'}'\left(\widetilde{\widetilde{A}}(x^{iv}_1,\ldots, x^{iv}_m)\right)= \widetilde{A}\left(i_{X'}'(x^{iv}_1),\ldots, i_{X'}'(x^{iv}_m)\right),$$
we continue with the same procedure, one variable at a time from right to left, as done in the proof of \cite[Lemma 1.2]{AroGalGarMae96}.

$(iii)\Rightarrow (iv)$ Recall that a functional in $X'$ is $w^*$-continuous if and only if it is $w^*$-continuous when restricted to bounded sets \cite[Theorem V.5.6]{dunford1988linear}.  Thus, if $\widetilde{A}$ is not $(w^*-w^*)$-continuous on the second variable, there exist $z_1,\ldots, z_m\in X''$, $e'\in E'$, $\varepsilon>0$ and a bounded net $(z_\alpha)\subseteq X''$ with $ z_\alpha \overset{w^*}{\to}z_2$   such that
\begin{equation}\label{ecu1}
\left| \widetilde{A}(z_1,z_\alpha,z_3,\ldots,z_m)(e')- \widetilde{A}(z_1,z_2,z_3,\ldots, z_m)(e')\right|>\varepsilon\qquad \forall\alpha.
\end{equation}

By the weak-star compactness of $B_{X^{iv}}$, take $(z_\beta)$ a subnet that is $\sigma(\cd{X},X''')$-convergent to some $\cd{x}\in\cd{X}$.  It is clear that 
\begin{equation}\label{ecu2}
i_{X'}'(\cd{x})=\sigma(X'',X')-\lim_\beta i_{X'}'(z_\beta) =\sigma(X'',X')-\lim_\beta z_\beta= z_2.
\end{equation}
By the definition of $\widetilde{\widetilde{A}}$,
$$
\widetilde{\widetilde{A}}(z_1,\cd{x},z_3,\ldots, z_m)= \sigma(X^{iv},X''')-\lim_\beta i_{E''} \circ \widetilde{A}(z_1,z_\beta,z_3,\ldots, z_m)
$$ then
\begin{eqnarray}
\label{ecu3}
i_{E'}'\circ\widetilde{\widetilde{A}}(z_1,\cd{x},z_3,\ldots, z_m)(e') &=& \widetilde{\widetilde{A}}(z_1,\cd{x},z_3,\ldots, z_m)(i_{E'}(e')) \\ \nonumber &=& \lim_\beta \widetilde{A}(z_1,z_\beta,z_3,\ldots, z_m)(e').
\end{eqnarray}

By \eqref{ecu1}, \eqref{ecu2} and \eqref{ecu3},
$$
i_{E'}'\circ\widetilde{\widetilde{A}}(z_1,\cd{x},z_3,\ldots, z_m)\neq \widetilde{A}(z_1,i_{X'}'(\cd{x}),z_3,\ldots, z_m)
$$ which contradicts $(iii)$.

$(iv)\Rightarrow (v)$ Given $z_1,\ldots, z_m\in X''$ we want to see that
$$\widetilde{A}(z_1,z_2,z_3,\ldots, z_m)= \widetilde{A}(z_2,z_1,z_3,\ldots, z_m).$$

For every $i$ take a net $(x_{\alpha_i})\subseteq X$ such that $ \wslim x_{\alpha_i}=z_i$. Then, by hypothesis
$$\widetilde{A}(z_1,z_2,z_3,\ldots, z_m)=   \wslim_{\,\,\,\,\,\,\,\,\,\,\,\,\,\,\,\,\,\,\,\,\,\,\alpha_2} \widetilde{A}(z_1,x_{\alpha_2},z_3,\ldots, z_m).$$
On the other hand, by definition of the  extension and the symmetry of $A$ we have
\begin{eqnarray*}
\widetilde{A}(z_2,z_1,z_3,\ldots, z_m)&=&   \wslim_{\alpha_2,\alpha_1,\alpha_3,\ldots, \alpha_m} A (x_{\alpha_2},x_{\alpha_1},x_{\alpha_3},\ldots, x_{\alpha_m})\\
&=&   \wslim_{\alpha_2,\alpha_1,\alpha_3,\ldots, \alpha_m} A(x_{\alpha_1},x_{\alpha_2},x_{\alpha_3},\ldots, x_{\alpha_m})\\
&=&   \wslim_{\,\,\,\,\,\,\,\,\,\,\,\,\,\,\,\,\,\,\,\,\,\,\alpha_2} \widetilde{A}(z_1,x_{\alpha_2},z_3,\ldots, z_m).\
\end{eqnarray*}

$(v)\Rightarrow(i)$ Since all permutations are compositions of  transpositions, it is enough to see that adjacent variables commute. That is, for any $k=1,\ldots,m-1$ and  $z_1,\ldots, z_m\in X''$ we have to show
$$\widetilde{A}(z_1,\ldots, z_k,z_{k+1},\ldots, z_m)=\widetilde{A}(z_1,\ldots, z_{k+1},z_{k},\ldots, z_m).$$

As before, for every $i$ take a net $(x_{\alpha_i})\subseteq X$ such that $  \wslim x_{\alpha_i}=z_i$. Then, by the definition of the  extension, the symmetry of $A$ and $(v)$ we have
\begin{eqnarray*}
\widetilde{A}(z_1,\ldots, z_k,z_{k+1},\ldots, z_m) &=&   \wslim_{ \alpha_1,\ldots,  \alpha_m} A(x_{\alpha_1},\ldots, x_{\alpha_k},x_{\alpha_{k+1}},\ldots, x_{\alpha_m})\\
&=&  \wslim_{ \alpha_1,\ldots, \alpha_m} A(x_{\alpha_{k}},  x_{\alpha_{k+1}}, x_{\alpha_1},\ldots, x_{\alpha_m})\\
&=&  \wslim_{ \alpha_1,\ldots, \alpha_{k-1}} \widetilde{A}(z_k, z_{k+1}, x_{\alpha_1},\ldots, x_{\alpha_{k-1}},z_{k+2},\ldots, z_m)\\
&=&  \wslim_{ \alpha_1,\ldots, \alpha_{k-1}} \widetilde{A}( z_{k+1}, z_k, x_{\alpha_1},\ldots, x_{\alpha_{k-1}},z_{k+2},\ldots, z_m)\
\end{eqnarray*}

On the other hand
\begin{eqnarray*}
\widetilde{A}(z_1,\ldots, z_{k+1},z_{k},\ldots, z_m) &=&   \wslim_{ \alpha_1,\ldots, \alpha_{k+1},\alpha_k,\ldots \alpha_m} A(x_{\alpha_1},\ldots, x_{\alpha_{k+1}},x_{\alpha_{k}},\ldots, x_{\alpha_m})\\
&=&  \wslim_{ \alpha_1,\ldots, \alpha_{k+1},\alpha_k,\ldots \alpha_m} A(x_{\alpha_{k+1}},  x_{\alpha_{k}}, x_{\alpha_1},\ldots, x_{\alpha_m})\\
&=&  \wslim_{ \alpha_1,\ldots, \alpha_{k-1}} \widetilde{A}(z_{k+1}, z_{k}, x_{\alpha_1},\ldots, x_{\alpha_{k-1}},z_{k+2},\ldots, z_m),
\end{eqnarray*}
which finishes the proof.
\end{proof}

Recall that, given a Banach space $X$, every $A\in \mathcal{L}(^mX;E)$ has a symmetric extension if and only if every bilinear map $B\in \mathcal{L}(^2X)$ has a symmetric extension. When we are working with a fixed  $A\in \mathcal{L}(^mX;E)$, a reformulation of $(v)$ in the previous result says that $A$ has a symmetric extension if and only if for any $z_3,\ldots, z_m\in X''$ the bilinear map $B\in \mathcal{L}(^2X; E)$ defined as $B(x_1,x_2)=\widetilde{A}(x_1,x_2,z_3,\ldots, z_m)$ has a symmetric extension.

As  mentioned above, on a symmetrically regular Banach space $X$  every symmetric bilinear  form $B\in \mathcal{L}(^2X)$ has a symmetric canonical extension $\widetilde{B}$. On non-symmetrically regular spaces this relationship (between regularity and symmetric extension) remains valid in the following sense: a symmetric bilinear  form $B\in \mathcal{L}(^2X)$ has symmetric  extension 
if and only if the linear operator $T\in \mathcal{L}(X,X')$ associated to $B$ is weakly compact (see the comments after  Lemma 8.1 in \cite{AroColGam91}). 

A symmetric $m$-linear form $A\in \mathcal{L}(^mX)$ if said to be \textit{regular} if its associated $(m-1)$-linear mapping $T_A\in \mathcal{L}(^{m-1}X,X')$ is weakly compact. So, a symmetric bilinear form is regular if and only if its canonical extension is symmetric. One implication of the previous sentence remains valid in the multilinear case: every regular symmetric $m$-linear form has symmetric  extension. Indeed, if $A\in \mathcal{L}_s(^mX)$ is regular then, by \cite[Theorem 1]{aron1997weakly},  there exist reflexive Banach spaces $Y_1, \dots Y_m$,  operators $T_i\in \mathcal{L}(X, Y_i)$ and  an $m$-linear mapping $D\in \mathcal{L}(Y_1\times \cdots \times Y_m)$  such that $A=D\circ (T_1\times \cdots\times T_m)$. Then, $\widetilde A= D\circ (T_1''\times \cdots\times  T_m'')$. Since $D$ is $w$-continuous on each variable and, by  the reflexivity of $Y_i$, each $T_i''$ is $(w^*-w)$-continuous, we get that $\widetilde A$ is $w^*$- continuous on each variable and then symmetric. 
However, the converse does not extend to symmetric $m$-linear forms, as the next example shows.

\begin{exa} There is a  non-regular trilinear mapping with symmetric canonical extension.

\rm Consider the space $X= c_0 \oplus_\infty \ell_2$ and the symmetric trilinear map $A\in \mathcal{L}(^3X)$ defined as
$$A((u,x),(v,y),(w,z))= \sum_{n=1}^\infty (x_ny_nw_n+x_nv_nz_n+u_ny_nz_n).$$

Since every operator from $X$ to $\ell_2$ is weakly compact (due to the reflexivity of $\ell_2$) and every operator from $X$ to $\ell_1$ is compact (due to Pitt's Theorem),  $X$ is a regular space and, hence, symmetrically regular. Then, the extension $\widetilde{A}$ is symmetric. Let us show that $A$ is not regular; that is, $T_A\in \mathcal{L}(^2X,X')$ given by
$$T_A((u,x),(v,y))= ((x_ny_n)_{n\in \zN}, (x_nv_n + u_ny_n)_{n\in \zN})$$
is not weakly compact. For this, just note  that
$$T_A((e_k,e_k), (e_k,e_k))=(e_k, 2\,e_k)$$
and that $\{(e_k, 2\,e_k)\}_{k\in \zN}\subset X'=\ell_1\oplus_1 \ell_2$ has no weakly convergent subnet, since $\{(e_k)\}_{k\in \zN}\subset \ell_1$ has no such a subnet.

\end{exa}

\bigskip

Notice that when $E$ is reflexive, Theorem \ref{th AB sim}  $(iii)$ tells us that
\begin{equation}\label{eq-representacion}\widetilde{\widetilde{A}}=\widetilde{A}\circ (i_{X'}'\times \cdots \times  i_{X'}').\end{equation}
In particular, this says that if $E$ is reflexive and $\widetilde{A}\in\mathcal L(^nX'',E)$ is symmetric then $\widetilde{\widetilde{A}}$ is also  symmetric. In  this case the third extension satisfies
\begin{equation}\label{eq-representacion2}
    \widetilde{\widetilde{\widetilde{A}}} = \widetilde{\widetilde{A}} \circ (i_{X'''}'\times \cdots \times  i_{X'''}')= \widetilde{A}\circ (i_{X'''}'\times \cdots \times  i_{X'''}') \circ (i_{X'}'\times \cdots \times  i_{X'}').
\end{equation}
That is, expression \eqref{eq-representacion} means that the second extension $\widetilde{\widetilde{A}}$  equals  the first extension composed with the natural projection from $X^{iv}$ onto $X''$, while \eqref{eq-representacion2} indicates that the third extension $\widetilde{\widetilde{\widetilde{A}}}$  coincides with the first extension composed with the natural projection from $X^{vi}$ onto $X''$.
Continuing like this, we obtain that the $k$-th extension $\widetilde{A}^{k}$ is equal to the first extension $\widetilde{A}$ composed with the natural projection from $X^{(2k)}$ onto $X''$. Consequently, if $E$ is reflexive and $\widetilde A$ is symmetric then $\widetilde A^k$ is symmetric, for all $k$. This  happens in particular in the scalar valued case. 

 In the language of the previous sections,  if $P\in\mathcal P(^mX)$ is the homogeneous polynomial associated to $A$ and  $\widetilde A$ is symmetric, the representation  \eqref{eq-representacion} implies that  $$\tilde{\tilde{\delta}}_{\cd{x} }(P)=\tilde{\delta}_{\pi(\cd{x})}(P)$$for every $\cd{x}\in \cd{X}$. For extensions to higher duals an analogous result holds. 
As a consequence, to distinguish evaluations at points of the bidual from evaluations at points of higher duals we need polynomials whose associated  symmetric multilinear mappings have non symmetric extension. 

\bigskip

Next, we show that in Theorem \ref{th AB sim} $(v)$ it is important that the variables that must commute are the first two. For that, we construct a trilinear form $A$ such that two variables of $\widetilde{A}$ commute, despite $\widetilde{A}$ not being symmetric. First, we prove a result about how to produce certain trilinear forms on $\ell_1$, which is a trilinear version of Lemma \ref{prop bilineal de ele1}.

\begin{lem}\label{lem trilineal de ele1} Given $a,b,c\in \ell_\infty$, there is a symmetric trilinear form
$A\in \mathcal L_s(^3\ell_1)$
such that, if $z_1\in\ell_1''$ and $z_2,z_3\in c_0^\bot\subseteq \ell_1''$, then
$$\widetilde{A}(z_1,z_2,z_3)=z_1(a)z_2(b)z_3(c).$$
\end{lem}
\begin{proof} The construction is similar to the one in Lemma \ref{prop bilineal de ele1} and we use some notation introduced at its proof. Define the trilinear form $\varphi\in \mathcal L(^3\ell_1)$ as 
\begin{eqnarray*}
\varphi(x,y,u)&=&\sum_{n= 1}^\infty \sum_{k=1}^n \sum_{l=1}^k c_nx_nb_ky_ka_lu_l= \sum_{n= 1}^\infty \sum_{k=1}^n c_nx_nb_ky_ka^{\leq k}(u) \\
&=&\sum_{n= 1}^\infty \sum_{l=1}^n\sum_{k=l}^n  c_nx_nb_ky_ka_lu_l =  \sum_{n= 1}^\infty \sum_{l=1}^n  c_nx_n a_lu_l (b^{\leq n} - b^{\leq l-1})(y) \\
&=&\sum_{k= 1}^\infty \sum_{l=1}^k\sum_{n=k}^\infty  c_nx_nb_ky_ka_lu_l =  \sum_{k= 1}^\infty \sum_{l=1}^k  b_ky_ka_lu_l c^{\geq k}(x). \
\end{eqnarray*}

Let us show that the symmetric trilinear form 
$$A(x,y,u):=\varphi(x,y,u)+\varphi(x,u,y)+\varphi(y,x,u)+\varphi(y,u,x)+\varphi(u,x,y)+\varphi(u,y,x)$$
meets the required conditions. Let $z_1\in\ell_1''$ and $z_2,z_3\in c_0^\bot\subseteq \ell_1''$. First observe that
$$z_3(a^{\leq k})=z_3(b^{\leq n} - b^{\leq l-1})=0  \,\,\,\text{ and }  z_3(c^{\geq k})=z_2(c).$$
Thus, if $(u^\gamma)\subseteq \ell_1$ is a net $w^*$-convergent to $z_3$, we have
\begin{eqnarray*}
\lim_{\gamma} \varphi (x,y,u^\gamma)&=& \sum_{n= 1}^\infty \sum_{k=1}^n c_nx_nb_ky_k z_3(a^{\leq k})=0\\
\lim_{\gamma} \varphi (x,u^\gamma, y)&=& \sum_{n= 1}^\infty \sum_{l=1}^n  c_nx_na_ly_l z_3(b^{\leq n} - b^{\leq l-1})=0\\
\lim_{\gamma} \varphi (u^\gamma,x, y)&=& \sum_{k= 1}^\infty \sum_{l=1}^k  b_kx_ky_lx_l z_3(c^{\geq k}) =z_3(c)\phi(x,y),\
\end{eqnarray*}
where $\phi(y,x)$ is the mapping defined in the proof of Lemma \ref{prop bilineal de ele1}.

Then, if $(x^\alpha), (y^\beta)\subseteq \ell_1$ are nets $w^*$-convergent to $z_1$ and $z_2$ respectively, and $B(z_1,z_2)$ is the bilinear form introduced in Lemma \ref{prop bilineal de ele1}, we have that
\begin{eqnarray*}
\widetilde{A}(z_1,z_2,z_3)&=& \lim_{\alpha, \beta, \gamma} A(x^\alpha, y^\beta, u^\gamma)\\
&=& \lim_{\alpha, \beta, \gamma} 
\Bigg(\varphi(x^\alpha, y^\beta, u^\gamma)+\varphi(x^\alpha, u^\gamma,  y^\beta)+\varphi( y^\beta,x^\alpha, u^\gamma)\\
& &+\varphi(y^\beta, u^\gamma, x^\alpha)+\varphi( u^\gamma, y^\beta, x^\alpha)+\varphi( u^\gamma, x^\alpha, y^\beta)\Bigg)\\
&=& \lim_{\alpha, \beta} 
\Bigg(0 + 0 +0 +0 +z_3(c)\phi(y^\beta, x^\alpha)+ z_3(c)\phi(x^\alpha,y^\beta,) \Bigg)\\
&=& z_3(c)\tilde{B}(z_1,z_2)=z_1(a)z_2(b)z_3(c)  \ \qedhere
\end{eqnarray*}  
\end{proof}

Following the same idea of the proof of Lemma \ref{lem trilineal de ele1}, it is not hard to show that given $a_1,\ldots, a_m\in \ell_\infty$, there is $A\in \mathcal{L}_s(^m\ell_1)$ 
such that, if  $z_2,\ldots, z_m\in c_0^\bot\subseteq \ell_1''$, then
$$\widetilde{A}(z_1, \ldots, z_m)=z_1(a_1)\cdots z_m(a_m).$$

Now, we show the promised example of a non-symmetric extension with two commuting variables.

\begin{exa} There exists a mapping $A\in\mathcal L_s(^3\ell_1)$ such that the second and third variable of $\widetilde{A}$ commute but $\widetilde{A}$ is not symmetric. 

\rm Take $A$ as in the Lemma \ref{lem trilineal de ele1} with $a,b\in \ell_\infty\setminus c_0$ linearly independent and $c=b$. Then it is clear that  $\widetilde{A}(z_1,z_2,z_3)\not=\widetilde{A}(z_2,z_1,z_3)$ for many $z_1,z_2,z_3\in\ell_1''$.

To see the commutativity of the last two variables note that by the construction of the canonical extension we have
$\widetilde{A}(z_1,z_2,z_3)=\widetilde{A}(z_1,z_3,z_2)$ when either $z_2$ or $z_3$ are elements in $\ell_1$. Now, since $\ell_1''=\ell_1\oplus c_0^\bot$, we only need to check that $z_2$ and $z_3$ commute in the case $z_2,z_3\in c_0^\bot$. This is immediate due to  Lemma \ref{lem trilineal de ele1}:
$$\widetilde{A}(z_1,z_2,z_3)=z_1(a)z_2(b)z_3(c)= z_1(a)z_2(c)z_3(b)=\widetilde{A}(z_1,z_3,z_2).$$
\end{exa}

\bigskip

Even when the canonical extension of a symmetric $m$-linear form $A$ is not symme\-tric, if one or more variables are taken in $X$, some permutations leave $\widetilde A$ unchanged. Let us take $m=3$ for simplicity. For $z_1\in X''$ and  $x_2, x_3\in X$,  $\widetilde A(z_1,x_2,x_3)$ coincides with $\widetilde A$ in any permutation of $\{z_1,x_2,x_3\}$. Also, for $z_1, z_2\in X''$ and  $x_3\in X$, it easily follows from the definition of the extension that \begin{equation}\label{eq-unaenX}
\widetilde A(z_1,z_2,x_3)=\widetilde A(z_1,x_3,z_2)=\widetilde A(x_3,z_1,z_2).   
\end{equation}We might be tempted to assert that, whenever a variable is in $X$, $\widetilde A$ is essentially symmetric.  However, the next example shows that if just one variable is in $X$, what we have in \eqref{eq-unaenX} is the best we can get. 
 
Recall that $\ell_1''=\ell_1\oplus c_0^\bot$.

\begin{exa} There is a mapping $A\in\mathcal L_s(^3\ell_1)$ such that $\widetilde{A}(\,\cdot\,,\,\cdot\,,z)$ is symmetric for any $z\in c_0^\bot$ but $\widetilde{A}(\,\cdot\,,\,\cdot\,,x)$ is not symmetric for many $x\in\ell_1$.
\end{exa}
\rm Take $A$ as in Lemma \ref{lem trilineal de ele1} with $a,c\in \ell_\infty\setminus c_0$ linearly independent and $b=a$. Arguing as in the previous example, it is easy to see that for any $z\in c_0^\bot$ we have that $\widetilde{A}(\,\cdot\,,\,\cdot\,,z)$ is symmetric. Now, let $x\in\ell_1$  such that $x(a)\neq 0$. In this case, if $z_1,z_2\in c_0^\bot$ satisfy $z_1(a)z_2(c)\neq z_2(a)z_1(c)$, then $\widetilde{A}(z_1,z_2,x)\neq  \widetilde{A}(z_2,z_1,x)$. Indeed, 
\begin{eqnarray*}
\widetilde{A}(z_1,z_2,x) &=& \widetilde{A}(z_1,x, z_2) 
=\widetilde{A}(x,z_1,z_2)
= x(a) z_1(a)z_2(c) \\
&\neq &    x(a) z_2(a)z_1(c)
= \widetilde{A}(x,z_2,z_1)= \widetilde{A}(z_2,z_1,x).
\end{eqnarray*}



\bibliography{main.bib}

\bibliographystyle{abbrv}

\end{document}